\documentclass[11pt, a4paper,leqno]{amsart}
\usepackage{amsmath,amsthm,amscd,amssymb,amsfonts, amsbsy}
\usepackage{latexsym}
\usepackage{txfonts}
\usepackage{exscale}
\usepackage{graphicx}
\usepackage{bbm}
\usepackage{romanbar}
\usepackage{cite}
\usepackage{enumerate}
\usepackage[utf8]{inputenc}

\usepackage[colorlinks,citecolor=red,pagebackref,hypertexnames=false]{hyperref}
\usepackage{color}

\calclayout
\allowdisplaybreaks

\theoremstyle{plain}
\newtheorem{theorem}[equation]{Theorem}
\newtheorem{lemma}[equation]{Lemma}

\theoremstyle{definition}
\newtheorem{definition}[equation]{Definition}

\theoremstyle{remark}

\numberwithin{equation}{section}

\newcommand{\RR}{{\mathbb{R}}}
\newcommand{\NN}{{\mathbb{N}}}

\newcommand{\eps}{\varepsilon}

\newcommand{\dint}{\int\!\!\!\int}

\newcommand{\dist}{\operatorname{dist}}

\newcommand{\re}{\mathbb{R}}
\newcommand{\rn}{\mathbb{R}^n}

\newcommand{\ree}{\mathbb{R}^{n+1}}

\newcommand{\N}{\mathbb{N}}

\newcommand{\dd}{\mathbb{D}}

\newcommand{\om}{\Omega}

\newcommand{\F}{\mathcal{F}}

\newcommand{\W}{\mathcal{W}}

\newcommand{\B}{\mathcal{B}}

\newcommand{\oo}{\mathcal{O}}

\newcommand{\mut}{\mathfrak{m}}

\newcommand{\pom}{\partial\Omega}

\newcommand{\tin}[1]{\scalebox{0.5}{#1}}
\newcommand{\X}{\textbf{X}}
\newcommand{\Y}{\textbf{Y}}
\newcommand{\Z}{\textbf{Z}}

\newcommand{\ulr}[1]{u\left(#1\right)}

\def\dashint{\Xint-}

\newcommand{\bestrel}[1]{\textnormal{\raisebox{-4pt}{$\stackrel{#1}{\sim}$}}}

\renewcommand{\emptyset}{\text{\textup{\O}}}

\DeclareMathOperator{\diam}{diam}

\def\div{\mathop{\operatorname{div}}\nolimits}

\newcommand{\vertiii}[1]{{\left\vert\kern-0.15ex\left\vert\kern-0.15ex\left\vert #1
		\right\vert\kern-0.15ex\right\vert\kern-0.15ex\right\vert}}

\def\Xint#1{\mathchoice
{\XXint\displaystyle\textstyle{#1}}%
{\XXint\textstyle\scriptstyle{#1}}%
{\XXint\scriptstyle\scriptscriptstyle{#1}}%
{\XXint\scriptscriptstyle%
\scriptscriptstyle{#1}}%
\!\int}
\def\XXint#1#2#3{{\setbox0=\hbox{$#1{#2#3}{%
\int}$ }
\vcenter{\hbox{$#2#3$ }}\kern-.6\wd0}}
\def\barint{\,\Xint -} 
\def\bariint{\barint_{} \kern-.4em \barint}
\def\bariiint{\bariint_{} \kern-.4em \barint}
\renewcommand{\iint}{\int_{}\kern-.34em \int} 
\renewcommand{\iiint}{\iint_{}\kern-.34em \int} 

\newcommand{\w}{\omega}

\begin{document}
\allowdisplaybreaks

\title[Characterization of Parabolic $L^p$ Solvability via Bounded Solution Estimates]{A Characterization of Solvability of the Parabolic $L^p$ Dirichlet Problem on Lipschitz Graph Domains Via Carleson Measure Estimates of Bounded Solutions}

\author{James Warta}
\address{
Department of Mathematics
\\
University of Missouri
\\
Columbia, MO 65211, USA}
\email{jw34r@umsystem.edu}

\author{Steve Hofmann}
\address{
Department of Mathematics
\\
University of Missouri
\\
Columbia, MO 65211, USA}
\email{hofmanns@umsystem.edu}

\thanks{}

\date{\today}


\begin{abstract}
In this paper, we show that if the bounded solutions to the parabolic Dirichlet problem on a Lipshitz-$\left[1,\frac{1}{2}\right]$ domain obey a Carleson measure estimate, then the corresponding parabolic measure on the boundary will belong to class $A^\infty$, which is equivalent to $L^p$ solvability for some $p<\infty$. This improves the existing literature which places additional assumptions on the parabolic uniform rectifiability or, equivalently, on the half-order time derivative of the function whose graph defines the boundary of the domain. 
\end{abstract}

\maketitle

\tableofcontents 

\section{Introduction}

 In this paper, we look to characterize the solvability of the Dirchlet problem with boundary data in $L^p$ for operators of the following form:
\begin{equation}\label{parabolic}Lu=\nabla\cdot A\nabla u-\partial_tu=0\end{equation}
where $A=A(X,t):\rn\times\re\rightarrow \re^{n^2}$ is a matrix-valued function that is uniformly elliptic, meaning that
\begin{equation}\label{ellip}
    \exists \lambda \in\RR^+: \quad \forall \xi\in\ree,\,\lambda^{-1}|\xi|^2\leq \xi^TA\xi\leq \lambda|\xi|^2
\end{equation}
where $\lambda$ is a positive constant, independent of the point $(X,t)$. The question at hand is to characterize the solvability of the following system of equations, known as the Dirichlet problem: \begin{align*}
    Lu&=0\quad\textnormal{in }\om\\
    u&\equiv f\quad \textnormal{on }\pom
\end{align*} The existence of such a solution encounters several obstructions, including the geometric regularity of $\pom$, and the analytic properties of both $A(X,t)$ and $f$, particularly when $f$ is not bounded. To this end, restrictions must be made, and in what follows they will be to restrict the class of domains to those laying above a Lipschitz graph. The characterization of solvability for boundary data $f\in L^p$ for some $p<\infty$ in such domains will be done by specifying a property of the solutions to the Dirichlet problem with boundary data $f$ that is continuous and compactly supported. The solutions to such a Dirichlet problem will be referred to as \textit{bounded solutions}. The condition on such bounded solutions $u$ will be:
\begin{equation}\label{CME}
    \sup_{\substack{\tin{\Y}\in\pom\\0<r<\diam\Omega}}\frac{1}{r^{n+1}}\dint_{I(\tin{\Y},r)\cap\Omega}|\nabla u |^2\delta(\X) \,d\X\leq C||u||_{L^{\infty}}^2
\end{equation}
where $C$ is a positive constant, independent of $u$, thus fixed across all bounded solutions. In the above $\delta(\X) $ denotes the distance from a point $\X$ to  the boundary $\pom$.  \par In what follows we characterize solvability by showing that if \eqref{CME} is true , then the corresponding parabolic measure $\w_L^{\X}$ with arbitrary pole $\X\in\Omega$ is of class $A_{\infty}$ with respect to the parabolic surface measure $\sigma$ on the boundary. By parabolic measure $\w_L^{\X}$ with pole at $\X\in\om$ we are referring to the unique measure such that for data $\phi$ on $\pom$ the function
\begin{equation}\label{soln}
    u(\X) =\dint_{\pom}\phi(\Y)d\w_L^{\X}(\Y)
\end{equation}
solves the corresponding Dirichlet problem of \eqref{parabolic} in $\om$ and with data $\phi$ on $\pom$. $A_\infty$ means that for any surface cube $\Delta=I(\X,r)\cap\pom$ and some $\beta\in(0,1)$ then there will exist an $\alpha\in(0,1)$ depending only on the given $\beta$ and global, fixed, constants such that for any surface subcube $\Delta'=I(\X',r')\cap\pom\subset\Delta$ and Borel subset $F\subset\Delta'$ we have the following implication: 
\begin{equation}\label{ainf}
    \frac{\w_L^{{\X}_\Delta^+}(F)}{\w_L^{{\X}_\Delta^+}(\Delta')}\leq \alpha \,\Rightarrow\,\frac{\sigma(F)}{\sigma(\Delta')}\leq \beta
\end{equation} 
where $\X_\Delta^+$ is the time-forward corkscrew point relative to $\Delta$, which will be discussed in the next section. \par
A large portion of the current literature surrounding the parabolic Dirichlet problem either relies on a parabolic notion of uniform rectifiability, or restricts to the special case of the heat equation, or $A\equiv I$. For a graph domain defined by a boundary function $\Psi$, this uniform rectifiability, or regularity, is equivalent to having $||\mathcal{D}^{\frac12}_t\Psi||_{BMO}<\infty$, where $\mathcal{D}^{\frac{1}{2}}_tf=(2\pi\tau\hat{f})\check{\;}$. Under this assumption, much has been proven; in fact the $A_\infty$ to Carleson measure estimate for bounded solutions in the regular case being equivalent to the $A_\infty$ condition was shown in \cite{DINDOS20171155} and a Carleson measure estimate on the coefficients, as opposed to the estimate on bounded solutions, implying solvability for data in $L^p$ was shown in \cite{dindoš2016}. For the specific case of the heat equation, the equivalence of boundary regularity, the $A_\infty$ property of the caloric measure, and $L^p$ solvability was shown in \cite{bortz2024}. Lastly, for general parabolic operators, regularity has been shown to imply a Carleson measure estimate for bounded solutions in \cite{URandCME}. Beyond the regularity case, less is known, but relevant to the discussion is \cite{Nystrom97}, which shows the equivalence of $L^p$ solvability and the $A_\infty$ property for general parabolic operators on Lipschitz domains and \cite{gens} which extends such a result to a class of 1-sided NTA domains. In this paper, we extend the conversation of the Carleson measure estimate condition and the $A_\infty$ condition to the realm of general parabolic operators and general Lipschitz graphs. In particular, we show that \eqref{CME} for bounded solutions is a complete characterization of the $A_\infty$ condition and, hence, $L^p$ solvability of the Dirichlet problem. \par To show that a Carleson measure estimate implies the $A_\infty$ property of the harmonic measure with respect to the parabolic surface measure, we look to the methods of \cite{Cavero_2020}, which are extensions of those originally found in \cite{kenig2014squarefunctionsainftyproperty}. This theory shows the equivalence of \eqref{CME} and \eqref{ainf} for elliptic generalizations of Laplace's equation in a more general class of domains, known as 1-sided chord arc domains, and we extend them to the parabolic setting on Lipschitz graph domains. To do this, we make use of parabolic analogues in \cite{Nystrom97} of elliptic estimates found in \cite{Cavero_2020}. We make note that while the estimates found in \cite{Nystrom97} are valid for domains which are parabolic Lipschitz cylinders, extending nicely to our setting, but they are only for solutions to parabolic equations in which $A$ is symmetric. To rectify this, we note that the proofs in \cite{Nystrom97} rely on the Harnack principle, energy estimates, Holder continuity, and the maximum principle of solutions to problems where $A(X,t)$ is always a symmetric matrix. All of these estimates have since been extended to general, non-symmetric parabolic problems.\par In the opposite direction, we use consequences of the $A_\infty$ property from \cite{URandCME} that are general enough to extend to our setting, and use them to provide estimates on the adjoint Green's function in relation to distance from the boundary of the domain. This allows for a an integration by parts over a region that contains the domain of integration to obtain an upper bound.
\section{Notations and Definitions}
In what follows, we characterize the space in $\mathbb{R}^{n+1}$ as $\RR^n\times\RR$, that is, a spatial plane $\rn$ crossed with a time dimension $\RR$. To this end, we denote $\X\in\ree$ as $\X=(X,t)=(x_0,x,t)$ with $X=(x_0,x)\in\rn$ as the spatial coordinate and $t$ as the time coordinate. We further divide the spatial plane $\rn$ into $\rn=\RR\times \RR^{n-1}$, denoting points $X=(x_0,x)$ as having a domain coordinate $x\in\RR^{n-1}$ and a so-called `graph coordinate' $x_0\in \RR$. For inequalities, if there exists a fixed, universal constant $C>0$, depending only on the structural properties of the problem, like dimension and geometry of the boundary, such that for two variable quantities $a$ and $b$ it's true that $a<Cb$ for some constant $C$ depending only on the so called `structural constants' of the problem, like the dimension or $\lambda$, then we write $a\,\bestrel{<}\,b$. If the constant depends on more quantities, say for example an arbitrary $\eps$, then we write $a\bestrel{<}_\eps b$ If $a\,\bestrel{<}\,b$ and $b\,\bestrel{<}\,a$, we denote this situation by writing $a\approx b$.
\par Now, make note of the homogeneity of the observed operator $L$: composing with the transformation $\chi:(X,t)\rightarrow(rX,r^2t)$ for some $r>0$ we obtain that:
$$L(u\circ\chi)=\div(A\nabla(\cdot))-\partial_t(\cdot)(u\circ\chi)=r\div(rA\mathbb{I}(\nabla u\circ\chi)-r^2(\partial_tu\circ\chi)= L(u)\circ\chi$$
This homogeneity is the motivation for the parabolic metric:
\begin{definition}[Parabolic Metric ($|||\cdot|||$)]For a given point $(X,t)\in\rn\times\RR$ we write: \begin{equation}\label{norm}|||(X,t)|||:=||X||+|t|^{\frac{1}{2}}\end{equation} 
where $||X||$ is the euclidean norm on $\ree$. 
\end{definition}
To ease notation, we use $\delta$ to denote a variety of distances under this metric. First, the distance between two points $\X$ and $\Y$ in $\ree$ under the parabolic metric will be denoted as $\delta(\X,\Y)$:
$$\delta(\X,\Y):=|||\X-\Y|||$$
the distance between a point $\X\in\ree$ and a set $E\subseteq\ree$ will be denoted as $\delta(\X,E)$:
$$\delta(\X,E):=\inf_{\Y\in E}\delta(\X,\Y)$$
and, as stated before, $\delta $ will be used to denote the distance from $\X$ to the boundary of the domain $\pom$, $\delta(\X) :=\delta(\X,\pom)$. 
 When we talk about a function $\Psi$ being Lipschitz, it is with respect to this parabolic metric. This property establishes the class of \textit{Lip-$[1,\frac{1}{2}]$ functions} and, correspondingly, Lipschitz graph domains with boundaries described by such functions are known as \textit{Lip-$[1,\frac{1}{2}]$ domains}:
 \begin{definition}[Lip-$\textnormal{[}1,\frac{1}{2}\textnormal{]}$ domains]  We say that $\om\in\ree$ is a \textit{Lip-$[1,\frac12]$ domain}, or a \textit{Lipschitz graph domain} if there exists a function $\Psi$ such that:
 $$\om=\left\{\X=(x_0,x,t)\in \ree: \Psi(x,t)>x_0\right\}$$
 where $||\Psi||_{Lip}<M=M(\pom)<\infty$. An implication of this is that $\pom$ is exactly the graph of the Lipschitz function $\Psi$.    
 \end{definition}
 Note that this class of domains is necessarily unbounded, but the results of this paper could be adapted to bounded domains via a standard localization argument. In the following arguments, we subdivide space in $\ree$ various ways, each with their own notation. First, a cube in $\ree$ centered at the point $\X=(X,t)$ with side length $2r$ will be denoted with either $I$ or $J$ with $I=I(\X,r)=I(X,t,r)$ defined as:
 $$I(\X,r):=\{\Y=(Y,s): |x_0-y_0|<3Mr,\,|y_i-x_i|<r
 \,\;\forall i\neq 0,\, |t-s|<r^2\}$$
 From this definition follows the definition of an arbitrary surface cube on $\pom$. We denote arbitrary surface cubes centered at $\X\in\pom$ with side length $2r$ as $\Delta(\X,r)$, defined as $\Delta(\X,r):=I(\X,r)\cap\pom$. Note that by the definition of $I(\X,r)$, $\Delta(\X,r)$ is a connected component. For both types of cubes, we denote the length of a cube to be half of its side length. Fro example, if $I=I(\X,r)$ then its length, denoted $l(I)$, would be $l(I)=r$. We will denote \textit{dyadic} surface cubes with either the letters $Q$ or $P$ and define them using the projection map to the space-time plane. Letting $\dd$ be the collection of dyadic cubes of $\re^{n-1}\times\re$ constructed using the parabolic metric, and define the map $\pi:(x_0,x,t)\rightarrow(x,t)$, the dyadic surface cubes $\dd(\pom)$ are the sets $Q\in\pom$ such that $\pi(Q)\in\dd$. This is to say, $\dd(\pom)$ are the sets in $\pom$ laying above the dyadic cubes in $\re^{n-1}\times \re$. For each $Q\in\dd(\pom)$ we define $\dd(Q)$ to be the set $\dd(Q):=\{Q'\in\dd(\pom): Q'\subset Q\}$. For each surface cube, we associate a set of three points in $\re^{n+1}$ known as the cubes \textit{corkscrew points}. 
\begin{definition}[Corkscrew points relative to $Q$,$\Delta$]{Given a cube $Q\in \dd(\pom)$, we define the \textit{time-forward corkscrew point relative to $Q=Q(y_0,y,s,r)$}}, $\X_Q^+$ to be the point: \begin{equation}\X_Q^+:=(y_0+k_1l(Q),\,y,\,s+k_2l(Q)^2)\end{equation}
where $k_1$ and $k_2$ are positive, structural constants to be chosen later and $l(Q)=2r$ is the length of $Q$. The \textit{time-backward} and \textit{centered corkscrew points relative to Q}, $\X_Q^-$ and $\X_Q$, respectively, are defined to be the points:
\begin{align*}
    \X_Q^-&:=(y_0+k_1l(Q),\,y,\,s-k_2l(Q)^2)\\
    \X_Q&:=(y_0+k_1l(Q),\,y,\,s)
\end{align*}
These points can also be defined for an arbitrary surface ball $\Delta=I(y_0,y,s,r)\cap\pom$, letting $l(\Delta)=l(I)$.
\end{definition}
These points are referenced in various inequalities and outside lemmas, and the constants $k_1$ and $k_2$ will be selected to simultaneously make use of them all. This can be done universally, and will only depend on the structural constants like $n$ and $||\Psi||_{Lip}$. \par The use of the parabolic norm necessitates the need for an appropriate analog of the Hausdorff measures of various dimension. We define the parabolic Hausdorff measure and state some of its basic properties here, drawing from \cite{BHHGNcorona}. A more complete discussion can be found in that text.
\begin{definition}[Parabolic Hausdorff Measure] Letting $\textnormal{diam}(E)$ denote the diameter of a set $E\subset\ree$ with respect to the parabolic metric the parabolic Hausdorff measure of dimension $\eta>0$ of $E$, $\mathcal{H}_p^\eta(E)$ is defined to be:
\begin{equation}
    \mathcal{H}_p^\eta(E)=\lim_{\delta\rightarrow 0^+}\mathcal{H}_{p,\delta}^\eta(E)
\end{equation}
where:
\begin{equation}
    \mathcal{H}_{p,\delta}^\eta(E):=\inf\left\{\sum_j\textnormal{diam}(E_j)^\eta\,:\,E\subset\bigcup_jE_j\,,\,\sup_j\textnormal{diam}(E_j)\leq\delta\right\}
\end{equation}
and:
\begin{equation}
    \textnormal{diam}(E):=\sup_{\X,\Y\in E}\delta(\X,\Y)
\end{equation}
\end{definition}
The first thing to make note of here is that, due to the structure of the parabolic metric, the homogeneous dimension of the entire space is $n+2$, even though the spacial dimension is $n+1$. Given that the boundary of the domain $\om$ is the graph of a Lipschitz function $\Psi$ over $\re^{n-1}\times\re$, which has homogeneous dimension $n+1$ we define the surface measure $\sigma$ to be: \begin{equation}\label{surf}
    \sigma:=\mathcal{H}^{n+1}_p|_{\pom}
\end{equation}
We end up working with this measure by analog; first we consider the following measure $\mu$:
\begin{equation}
    \label{mu}\mu(E):=\int_\re\int_{\re^n\times\{t\}}\mathbf{1}_{E}(X,t)d\mathcal{H}^{n-1}d\mathcal{H}^1(t)
\end{equation}
and let $\sigma^s=\mu|_{\pom}$. It is a result of \cite{BHHGNcorona} that in the case of a Lip-$[1,\frac12]$ graph domain $\sigma^s\approx\sigma$. This is due to a more primary result, and the fact that a Lip-$[1,\frac12]$ graph is \textit{parabolic Ahlfors-David regular} with respect to $\sigma^s$. 
\begin{definition}[Parabolic Ahlfors-David regular] We say that a closed set $E\subset \ree$ is a parabolic Ahlfors-David regular (ADR) set with respect to a measure $\mut$ if for all $r>0$, $\X\in E$ that: \begin{equation}\label{adr}\mut(B(\X,r)\cap E)\approx r^{n+1}\end{equation}
\end{definition}
Note that ADR with respect to balls extends to ADR w/respect to cubes, that is, it is equivalent to replace $B(\X,r)$ with $I(\X,r)$ in \ref{adr}. Lastly, we make use of a construction seen in \cite{Cavero_2020} known as a \textit{Good $\epsilon_0$-cover}:
\begin{definition}[Good $\epsilon_0$-cover] Let $E\subset\ree$ be an $n$-dimensional ADR set. Fix $Q_0\in\dd(E)$ and let $\mu$ be a regular Borel measure on $Q_0$. Given $\epsilon_0\in (0,1)$ and a Borel set $F\subset Q_0$, we say a \textit{good $\epsilon_0$-cover from F to $Q_0$ with respect to $\mu$} is a collection $\{\oo_l\}_{l=1}^k$ of Borel subsets of $Q_0$ with pairwise disjoint families $\F_l=\{Q_i^l\}\subset \dd(Q_0)$ such that:
\begin{itemize}
    \item[(a)] $F\subset \oo_k\subset\oo_{k-1}\subset\cdots\subset\oo_2\subset\oo_1\subset Q_0$ 
    \item[(b)] $\oo=\bigcup_{Q_i^l\in \F_l}Q_i^l,\quad\quad1\leq l\leq k$
    \item[(c)] $\mu(\oo_l\cap Q_i^{l-1})\leq \eps_0\mu(Q_i^{l-1})\quad\quad\forall Q_i^{l-1}\in\F_l,\quad2\leq l\leq k$
\end{itemize}
Where $k\in\N$ is a fixed integer, which we will refer to as the length of the cover in question.
\end{definition}
\section{Outside Lemmas}
The first set of outside lemmas to discuss surround the good $\eps_0$-covers defined at the end of the previous section. The proofs of these rely only on the structure of $\rn$ and the fact that a given measure $\mu$ is Borel, the reader is referred to \cite{Cavero_2020} for the proofs. The first surrounds their existence:
\begin{lemma}
    
Let $E\subset\ree$ be a ADR set with respect to the parabolic surface measure and have parabolic Hausdorff dimension $n+1$ and fix $Q_0\in\dd(E)$. Let $\mu$ be a regular Borel measure on $Q_0$ and assume further that it is dyadically-doubling on $Q_0$. For every $0<\eps_0<e^{-1}$, if $F\subset Q_0$ is a Borel subset with $\mu(F)<\alpha \mu(Q_0)$ and $0<\alpha<\frac{\eps_0^2}{2C_\mu^2}$ then $F$ has a good $\eps_0$-cover with respect to $\mu$ of length $k_0=(\alpha,\eps_0)\in\N$ where $k_0\geq 2$ and $k_0\approx \frac{\ln(\alpha)}{\ln(\eps_0)}$
\end{lemma} 
\begin{lemma}If $\{\oo_l\cap Q_i^m\}_{l=1}^k$ is a good $\eps_0$-cover of $F$ with respect to $\mu$ of length $k\in\NN$ then $$\mu(\oo_l\cap Q_i^m)\leq \eps_0^{l-m}\mu(Q_i^m),\quad\quad \forall Q_i^m\in \F_m,\quad\quad 1\leq m\leq l\leq k$$
\end{lemma}

The next set of lemmas are quantitative estimates for parabolic solutions and are restatements of those found in \cite{Nystrom97} to fit within the definition of cubes found here. The first such estimate is a parabolic Harnack inequality for solutions to parabolic operators on Lipschitz domains:
\begin{lemma}[Harnack Inequality]\label{Harnack}For $\X=(X,t)\in\pom$, $\om$ a Lip-[$\it{1,\frac12}$] domain, and $a,r>0$ let
$$J_a(\X,r):=\left\{\Y=(Y,s)\in\om\,:\,\delta(\Y)>ar,\,\,\delta(\X,\Y)<\left(1+\frac
1a\right)r,\,s>t+(1+a)r^2\right\}$$
Then, if $Lu=0$ in $\om$ and $u=0$ on $\pom\setminus Q(\X,r)$, there exists $C_a>0$ depending only on structural constants and $a$ such that: 
$$\sup_{J_a(\X,r)}u\leq C_a\inf_{J_a(\X,r)} u$$
\end{lemma}
We also draw on a doubling estimate of the parabolic measure on $\pom$:
\begin{lemma}[Doubling of the Parabolic Measure]\label{doubling} For $a>0$, $\X=(X,t)\in\pom$, with $\om$ a Lip-[$\it{1,\frac12}$] domain, let $$\Xi_a(\X) :=\left\{(Y,s)\in\om|\,||X-Y||^2<s-t\right\}$$
Then there exists a constant $\beta=\beta(M)$ depending only on the Lipschitz constant of the boundary such that if $\Y=(Y,s)\in\Gamma_a(\X) $ with $s>t+\beta r^2$ for given $r>0$ then:
$$\w_L^{\Y}(Q(\X,2r))\,\bestrel{<}\,\w_L^{\Y}(Q(\X,r))$$
For all $\Y\in\Gamma_a(\X)$ with the implicit constants of the inequality depending only on $a$ and the structural constants.
\end{lemma}
The last outside lemma we pull from \cite{Nystrom97} is the pole-change formula:
\begin{lemma}[Pole Change Formula]\label{polechange}
    With the same hypotheses as lemma 3.4, and any Borel set $E\subset Q=Q(\X,r)\in\dd(\pom)$ we have that: $$\frac{\w_L^{\Y}(E)}{\w_L^{\Y}(Q)}\approx \w_L^{\X_Q^+}(E)$$
    For all $\Y\in\Gamma_a(\X)$ with the implicit constants of the inequality depending only on $a$ and the structural constants.
\end{lemma}
The last outside lemma we make use of is a Bourgain estimate of the Caloric measure from \cite{gens}:
\begin{lemma}[Bourgain-Type Estimate]\label{bourg}
    There exist constants $\gamma,\eta>0$ such that for all $r>0$, $\X_0\in\pom$, if $\Y\in \Delta=\Delta(\X_0,\gamma r)$ then:
    $$\w_L^{\Y}\left(\Delta\right)\geq \eta$$
\end{lemma}
\section{CME Implies $A_\infty$}
In this section we prove the direction:
\begin{theorem}[CME implies $A_\infty$]\label{mainthm}
    Let $\om\in\ree$ be a a Lipschitz graph domain, and let $L$ be a uniformly-parabolic operator. If every bounded solution $u$ to the $L^p$ Dirichlet problem corresponding to $L$ satisfies the following inequality:
    \begin{equation}\label{4cme}
        \sup_{\substack{\Y\in\pom\\0<r<\diam\Omega}}\frac{1}{r^{n+1}}\dint_{I(\Y,r)\cap\Omega}|\nabla u |^2\delta(\X) \,d\X\leq C
    \end{equation}
    for some $C\in\re$, independent of the given bounded, solution $u$, then for every for every $\beta\in(0,1)$ there will exist an $\alpha\in(0,1)$ such that for every surface cube $\Delta_0$, surface subcube $\Delta\subset\Delta_0$ and Borel set $F\subset \Delta$ the following implication holds:
    \begin{equation}\label{4a}\frac{\w_L^{\X_{\Delta_0}^+}(F)}{\w_L^{\X_{\Delta_0}^+}(\Delta)}\leq \alpha \,\Rightarrow\,\frac{\sigma(F)}{\sigma(\Delta')}\leq \beta
    \end{equation}
\end{theorem}
To prove this, we first show that it will suffice to prove a stronger statement on the dyadic cubes:
\begin{lemma}\label{oncube}
    Assuming there exists a $C\in\re$ such that $\eqref{4cme}$ holds independent of any bounded solution $u$ of $\eqref{parabolic}$, then for any $\beta\in(0,1)$ there exists an $\alpha\in(0,1)$ such that given $Q_0\in\dd(\pom)$, for any $Q_0\in\dd(\pom)$ and every Borel subset $F\subset Q_0$ where: $$\frac{\w_L^{\X_{Q_0}^+}(F)}{\w_L^{\X_{Q_0}^+}(Q_0)}\leq \alpha$$ it will always be the case that: $$\frac{\sigma(F)}{\sigma(Q_0)}\leq \beta$$
\end{lemma}
Assuming this lemma for now, we provide a proof of \ref{mainthm}:
\begin{proof}[Proof of \ref{mainthm}] We first make use \ref{oncube} to formulate a dyadic $A_\infty$ condition. Note the following estimate: via the pole change formula \ref{polechange}, Harnack inequality \ref{Harnack} and Bourgain estimate \ref{bourg} we see that for a fixed $Q_0,Q^0\in\dd(\pom)$ and Borel subset $F$ such that $F\subset Q_0\subset Q^0$ we have for some constant $C$ that: \begin{equation}\label{aoncube}
    \frac{1}{C}\frac{\w_L^{\X_{Q_0}^+}(F)}{\w_L^{\X_{Q_0}^+}(Q_0)}\leq\frac{\w_L^{\X_{Q^0}^+}(F)}{\w_L^{\X_{Q
    ^0}^+}(Q_0)}
\end{equation}
Now we proceed to prove \ref{mainthm} by contraposition and note that if $\frac{\w_L^{X_{Q^0}^+}(F)}{\w_L^{X_{Q^0}^+}(Q_0)}<\frac{\alpha}{C}$ then \ref{oncube} and \ref{aoncube} yield that $\frac{\sigma(F)}{\sigma(Q_0)}<\beta$. With this, take any $\beta\in(0,1)$ and any surface ball $\Delta_0=\Delta(\X_0,r_0)$. Then take any sub-surface cube $\Delta=\Delta(\X,r)\subset\Delta_0$, and Borel subset $F\subset \Delta$ such that $\sigma(F)>\beta\sigma(\Delta)$. By Ahlfors regularity, it's possible to select a constant $K>1$ to create a family of cubes $\F_K:=\left\{Q\in\dd(\pom):Q\cap\Delta\neq\emptyset,\frac{r}{4K}<l(Q)<\frac{r}{2K}\right\}$ so that $\Delta\subset\bigcup_{\F_K}Q\subset 2\Delta$. In light of the fact that $\sigma(F)>\beta\sigma(\Delta)$, the pidgeon-hole principle and Ahlfors regularity give the existence of a$Q_0\in\F_K$ and a constant $K'>1$ for which $\frac{\sigma(F\cap Q_0)}{\sigma(Q_0)}>\frac{\beta}{K'}$. Thus, by contrapositive of \ref{oncube}, we have that $\frac{\w^{X^+_{Q_0}}(F\cap Q_0)}{\w_L^{X^+_{Q_0}}(Q_0)}>\alpha$ for some $\alpha\in(0,1)$.Also observe that by the doubling property of $\w_L^{\X_{\Delta_0}}(\cdot)$, we have that $\w_L^{\X_{\Delta_0}^+}(\Delta)\approx \w_L^{\X_{\Delta_0}^+}(Q_0)$. Now, let $Q^0\in\dd(\pom)$ be the unique dyadic cube containing $Q_0$ such that $\frac
{r_0}{2}<l(Q^0)\leq r_0$ so that, by definition of corkscrew points, $\delta(\X_{Q^0})\approx\delta(\X_{\Delta_0})$. Combining these facts together we have the following:
\begin{equation}
    \frac{\w_L^{\X_{\Delta_0}}(F)}{\w_L^{\X^+_{\Delta_0}}(\Delta)}>\frac{\w_L^{\X^+_{\Delta_0}}(F\cap Q_0)}{\w_L^{\X^+_{\Delta_0}}(\Delta)}\approx \frac{\w_L^{\X^+_{\Delta_0}}(F\cap Q_0)}{\w_L^{\X^+_{\Delta_0}}(Q_0)}\approx \frac{\w_L^{\X^+_{Q^0}}(F\cap Q_0)}{\w_L^{\X^+_{Q^0}}(Q_0)}>\frac{\w_L^{\X^+_{Q_0}}(F\cap Q_0)}{\w_L^{\X^+_{Q_0}}(Q_0)}\frac{1}{C}>\frac{\alpha}{K''}=\tilde{\alpha}
\end{equation}
Thus for any $\beta\in(0,1)$ and $\Delta_0$ there exists an $\tilde{\alpha}\in(0,1)$ such that if $\frac{\sigma(F)}{\sigma(\Delta)}>\beta$ then $\frac{\w_L^{X^+_{\Delta_0}}(F)}{\w_L^{X^+_{\Delta_0}}(\Delta)}$ for arbitrary $F\subset\Delta\subset\Delta_0$
\end{proof}
The proof of \ref{oncube} relies on another lemma, whose statement requires some definitions. First let $A(\X,r,R)$ be the cubic annulus centered at $\X$; $A(\X,r,R):=I(\X,R)\setminus I(\X,r)$. Further, for $Q\in\dd(\pom)$, we define the cube $\tilde{Q}\in\dd(Q)$ to be the dyadic descendant of $Q$ that contains the center of $Q$ and has side length $\eta l(Q)$ for some $\eta=2^{-k}<1$ to be chosen later. In addition, given some $\Y\in Q$ let $P_{{\Y}}\in\dd(Q)$ be the descendant with side length $\eta l(Q)$ but containing the point $\Y\in Q$. Now, for any cube $Q\in \mathbb{D}(\partial\Omega)$, $\Y\in Q$, and $\eta\in(0,1)$ define the following region:
$$(\Lambda^\eta)^\star_Q(\Y):=\{\X\in\Omega:\delta(\X) K_1\geq ||\X-\Y||\}\cap A(\Y,\eta^3l(Q),K_2l(Q))$$
and note that, if the constants $K_1$ and $K_2$ are chosen appropriately with respect to $M$, $k_1$ and $k_2$, this is a connected Whitney region which contains the corkscrew points of ${\Tilde{Q}}$ and ${\Tilde{P}_{{\Y}}}$, where $\Tilde{P}_{{\Y}}$ is to $P_{{\Y}}$ as $\Tilde{Q}$ is to $Q$. Furthermore, define the Lipschitz cone $\Gamma^\star_{Q}(\Y)$ to be: $$\Gamma^\star_Q(\Y):=\{\X\in\Omega: K_1\delta(\X) \geq||\X-\Y||\}\cap B(\Y,K_2l(Q))$$
Now, for some $\eta$, $Q$, and $y\in Q$ let $W^\eta_Q(y)$ be the collection of cubes in the Whitney decomposition of $\Omega$ that meet $(\Lambda^\eta)^\star_Q(y)$ and let $I^\sigma:=(1+\sigma)I$ be a concentric fattening of a given cube $I$. Then, there exists a $\sigma>0$ depending only on structural constants such that the region: $$\Lambda^\eta_Q(\Y)=\bigcup_{I\in W^\eta_Q(\Y)}I^\sigma$$
is entirely contained in $\Omega$ for any $\eta$, $Q$, and $\Y\in Q$. Define $\Gamma_Q(\Y)$ in relation to $\Gamma^\star_Q(\Y)$ in the same fashion
\\
\begin{lemma}\label{Sqfcnest}
    There exists constants $0<\eta<1$, $0<\alpha_0<1$, depending only on the Lipschitz constant of the boundary and the ellipticity of L, and a $C_\eta$ depending further on the $\eta$ chosen to define $\Lambda_Q^\eta$ such that for every $0<\alpha<\alpha_0$, for every $Q_0\in\mathbb{D}(\partial\Omega)$ and Borel set $F\subset Q_0$ satisfying $w_L^{\X_{Q_0}^+}(F)\leq \alpha w_L^{\X_{Q_0}^+}(Q_0)$ there exists a Borel set $S\subset Q_0$ such that the bounded weak solution $u\left(\X\right)=w_L^{\X}(S)$ satisfies:
    $$S_{Q_0}u\left(\X\right):=\left(\int_{\Gamma_{Q_0}\left(\tin{\X}\right)}|\nabla u(\Y)|^2\delta(\Y)^{1-n}d\Y\right)^{\frac{1}{2}}\geq C_\eta^{-1}(\ln(\alpha^{-1}))^{\frac{1}{2}}$$ for every $\X\in F$
\end{lemma}
In the following proof of \ref{Sqfcnest}, we make use of the following lemma, to be proven subsequently. The definition of $\tilde{Q}_i^l$ and $\tilde{P}_i^l$ will be given in the proof of \ref{Sqfcnest}.
\begin{lemma}\label{lowerS}
     Defining $u\left(\X\right)$ as above, for sufficiently small $\eta$ chosen to define $\Lambda_Q^\eta$ and $\epsilon_0$ chosen to define the good $\eps_0$-cover, both depending only on structural constants, it is true that for each $\Y\in F$, $1< l\leq k-1$, $$\left|u\left(\X_{\tilde{Q}_i^l}^+\right)-u\left(\X_{\tilde{P}_i^l}^+\right)\right|\geq \frac{c_0}{2}$$ where $c_0$ is a constant depending on the Ahlfors regularity property of $\Omega$ and the ellipticity of $L$
\end{lemma}
\begin{proof}[Proof of Lemma \ref{Sqfcnest}] Let $\eta$ be small enough so that Lemma 4.8 holds. For a given $Q_0\in\mathbb{D}(\partial\Omega)$ note that $\w=\w_L^{\X_{Q_0}^+}$ is a dyadically-doubling, (with constant $C_0$ depending only on structural constants) regular, Borel measure on $\partial\Omega$. Let $0<\epsilon_0<1$ and $0<\alpha<\epsilon_0^2/(2C_0)^2$ to be explicitly chosen later, and let $F$ be the given Borel subset of $Q_0$. By properties of good-$\epsilon_0$ covers, there exists a good-$\epsilon_0$ cover from $F$ to $Q_0$ with respect to $\w$ of length $k\approx \frac{\ln(\alpha)}{\ln(\epsilon_0)}$. Let $O_l$ be the collection of Borel sets such that $F\subset O_k\subset \cdots\subset O_1\subset Q_0$ where $O_l=\bigcup_{Q_i^l\in\mathcal{F}_l}Q_i^l$. From here we construct the set $S$ by letting $\Tilde{O}_l:=\bigcup_{Q_i^l\in\mathcal{F}_l}\Tilde{Q}_i^l$ and defining $S:=\bigcup_{j=2}^k(\Tilde{O}_{j-1}\setminus O_j)$. Now, noting that this union is disjoint, let $u(\X):= \w_L^\X(S)$ and observe that: 
\begin{equation}    
u\left(\X\right)=\int_{\partial \Omega}\mathbf{1}_S(\Y)d\w_L^{\X}(\Y)=\sum_{j=2}^k\w_L^{\X}(\Tilde{O}_{j-1}\setminus O_j)
\end{equation}
Now, for each $\Y\in F$ and $1\leq l\leq k$ there exists a unique $Q_i^l=Q_i^l(\Y)$ containing $\Y$ and let $P_i^l$ be the corresponding $P_{\tin{\Y}}$ containing $\Y$ to this $Q_i^l(\Y)$ as described above. Now, by construction, both $\X_{\Tilde{Q}_i^l}^+$ and $\X_{\Tilde{P}_i^l}^+$ are contained in $\Lambda^\eta_{Q_i^l}(\Y)$. In observation of the fact that this is a Whitney region, we let $I_{\Tilde{Q}_i^l}$ and $I_{\Tilde{P}_i^l}$ be the Whitney cubes that contain $\X_{\Tilde{Q}_i^l}^+$ and $\X_{\Tilde{P}_i^l}^+$, respectively. Using \ref{lowerS}, the fact that parabolic solutions have enough regularity to obey an interior Poincare estimate \cite{Brown_1990}, Moser local boundedness estimates \cite{Moser}, and the fact that the inclusion $\X\in \Lambda^\eta_{Q_i^l}(\Y)$ implies $\delta(\X) \approx_\eta l(Q_i^l)$, with the inequality constants additionally depending on $\eta$, we obtain that:
\begin{align}
\frac{c_0}{2}&\leq \left|\ulr{\X_{\Tilde{P}_i^l}^+}-u_{\Lambda^\eta_{Q_i^l}(\tin{\Y})}\right|+\left|\ulr{\X_{\Tilde{Q}_i^l}^+}-u_{\Lambda^\eta_{Q_i^l}(\tin{\Y})}\right|\\
&\bestrel{<}\left(\dashint_{I^\sigma_{\Tilde{Q}_i^l}}\left|\ulr{\X}-u_{\Lambda^\eta_{Q_i^l}(\tin{\Y})}\right|^2d\X\right)^{\frac{1}{2}}+\left(\dashint_{I^\sigma_{\Tilde{P}_i^l}}\left|\ulr{\X}-u_{\Lambda^\eta_{Q_i^l}(\tin{\Y})}\right|^2d\X\right)^{\frac{1}{2}}\\
&\leq C_\eta \left(l(Q_i^l)^{-(n+1)}\int_{\Lambda^\eta_{Q_i^l}(\tin{Y})}\left|\ulr{\X}-u_{\Lambda^\eta_{Q_i^l}(\tin{Y})}\right|^2 d\X\right)^{\frac{1}{2}}\\
&\leq C_\eta\left(\int_{\Lambda^\eta_{Q_i^l}(\tin{Y})}\left|\nabla \ulr{\X}\right|^2 \delta\left(\X\right)^{1-n}d\X\right)^{\frac{1}{2}}
\end{align}
Now, noting that for each $l$, $Q_i^l(\Y)$ is a dyadic cube that contains $
\Y$, thus the collection $\{Q_i^l(\Y)\}_l$ is a family of non-trivially intersecting dyadic cubes, meaning it is a nested family of cubes. By properties of the good-$\epsilon_0$ cover, we have that for any valid $l$ that $\w(Q_i^l\cap Q_i^{l-1})\leq \w(Q_i^l\cap O^{l-1})\leq \epsilon_0 w(Q_i^l)$. By the doubling properties of $\w$, we can specify a small enough $\epsilon_0$ such that $Q_i^l$ and $Q_i^{l-1}$ are of different length scales and $\Lambda^\eta_{Q_i^l}(\Y)$ and $\Lambda^\eta_{Q_i^{l-1}}(\Y)$ distinct Whitney regions, thus the $\Lambda^\eta_{Q_i^l}(\Y)$ have bounded overlap in $\Gamma_{Q_0}$. From this bounded overlap, we obtain that:
\begin{equation}\frac{c^2_0}{4}\frac{\ln(\alpha)}{\ln(\epsilon_0)}\simeq \frac{c^2_0}{4}(k-1)\leq C_\eta \sum_{l=1}^{k-1}\int_{\Lambda^\eta_{Q_i^l}(\tin{\Y})}\left|\nabla u \right|^2 \delta(\X) ^{1-n}d\X\leq C_\eta (S_{Q_0}u(\Y))^2 \cdot\end{equation}
\end{proof}
Now the only proof left before the proof of \ref{oncube} is the proof of \ref{lowerS}.
\begin{proof}[Proof of \ref{lowerS}]
 First, we need to estimate $u(\X_{\tilde{Q}_i^l}^+)$. Noting that $\delta(\X_{\Tilde{Q}_i^l}^+,Q_i^l)\bestrel{<}\eta l(Q_i^l)$, by Holder continuity at the boundary we have for some $\gamma$ and constant $C$ that:
\begin{align}
    u(\X_{\Tilde{Q}_i^l}^+)\leq \w_L^{\X_{\Tilde{Q}_i^l}^+}(\partial\Omega\setminus Q_i^l)+\w_L^{\X_{\Tilde{Q}_i^l}^+}(S\cap Q_i^l)\leq C\eta^{\gamma}+\w_L^{\X_{\Tilde{Q}_i^l}^+}(S\cap Q_i^l)
\end{align}
Following the construction of the good $\epsilon_0$-cover, write: 
\begin{equation}\w_L^{\X_{\Tilde{Q}_i^l}^+}(S\cap Q_i^l)=\sum_{j=l+2}^k\w_L^{\X_{\Tilde{Q}_i^l}^+}(Q_i^l\cap(\Tilde{\mathcal{O}}_{j-1}\setminus\mathcal{O}_j))+w_L^{\X_{\Tilde{Q}_i^l}^+}(Q_i^l\cap(\Tilde{\mathcal{O}}_{l}\setminus\mathcal{O}_{l+1}))\end{equation}
To estimate from above further, we'll first need to pass from $\w_L^{\X_{\Tilde{Q}_i^l}^+}$ to $\w_L^{\X_{Q_i^l}^+}$ via the Harnack inequality.
We claim that if $\Z_{Q_i^l}$ is the center of $Q_i^l$, with $r=l(Q_i^l)$ and $a=\frac{\eta}{\sqrt{2(k_1^2+k_2)}-1}$ then, $\X_{\Tilde{Q}_i^l}^+,\X_{Q_i^l}^+\in J_a(\Z_{Q_i^l},r)$. To see this, first note that $\delta(\X_{\Tilde{Q}_i^l}^+)<\delta(\X_{Q_i^l}^+)$ so it will suffice to show that $\delta(\X_{\Tilde{Q}_i^l}^+)>ar$. To do this, write $\Z_{Q_i^l}=(y_0,y,s)$ and define:
$$V_{\tin{\Z}_{Q_i^l}}:=\left\{(y_0+M|(x,t)-(y,s)|,x,t):(x,t)\in\mathbb{R}^{n}\right\}$$ and note that $\delta(\X_{\Tilde{Q}_i^l}^+,\partial\Omega)\geq \delta(\X_{\Tilde{Q}_i^l}^+,V_{\tin{\Z}_{Q_i^l}})$ further, for any $(x,t)\in\mathbb{R}^{n}$, we have that:
\footnotesize
\begin{align}
    |(\overline{y}+k_1l(\Tilde{Q}_i^l),y,s+k_2l(\Tilde{Q}_i^l)^2)-(\overline{y}+M|(y,s)-(x,t)|,x,t)|^2&\geq|k_1l-M|(y,s)-(x,t)||^2+|y-x|^2+|s-t|+k_2l^2\\
    &=:(k_1l-Mz)^2+z^2+k_2l^2
\end{align}
\normalsize
Which is minimized for $z=\left(\frac{Mk_1}{m^2+1}\right)l$ with value $l\sqrt{\frac{k_1^2}{M^2+1}+k_2}$. Now, note that $\delta(\X_{\Tilde{Q}_i^l}^+,\Z_{Q_i^l})\leq\delta(\X^+_{Q_i^l},\Z_{Q_i^l})$, so we only need to find sufficient conditions for $a$ on the latter term. Observe that:
\begin{align}
    |\delta(\X_{Q_i^l},\Z_{Q_i^l})|^2&=\left|\left|\left|(k_1l(Q_i^l),0,k_2l(Q_i^l))\right|\right|\right|^2\\
    &\leq2(k_1^2l(Q_i^l)^2+k_2l(Q_i^l)^2)=2l(Q_i^l)^2(k_1^2+k_2)
\end{align}
thus, if we need \begin{equation}l(Q_i^l)\sqrt{2(k_1^2+k_2)^2}\leq \left(1+\frac{1}{a}\right)l(Q_i^l)\end{equation} it follows that $a=\eta\frac{1}{\sqrt{2(k_1^2+k_2)}-1}$ suffices for both conditions simultaneously.\\
\\
Choosing such an $a$, which depends on $\eta$, means that the application of \ref{Harnack} will incur a constant depending on $\eta$ and we have that:
\begin{equation}\w_L^{\X_{\Tilde{Q}_i^l}^+}(Q_i^l\cap(\Tilde{O}_{j-1}\setminus O_j))\leq C_\eta \w_L^{\X_{Q_i^l}^+}(Q_i^l\cap(\Tilde{O}_{j-1}\setminus O_j))\end{equation}
Now, we seek to pass from $\w_L^{\X_{Q_i^l}^+}$ to $\w_L^{\X_{Q_0}^+}$ to apply the properties of the good $\epsilon_0$ cover. To do this we apply the pole change formula. To do that we show $\X_{Q_0}^+\in \Xi_a(\Z_{Q_i^l})$ for some $a>0$. To this end, we will write $\X_{Q_0}^+$ as $(X,\tau)=(x_0+k_1l(Q_0),x,t+k_2l(Q_0)^2)$, $(Y,s)=(y_0,y,s)=\Z_{Q_i^l}$ as the center of $Q_i^l$, and $I_r=Q_i^l$. First observe that:
\begin{align}
      |X-Y|^2&=|(x_0+k_1l(Q_0),x)-(y_0,y)|^2\\
    &=|x_0-y_0+k_1l(Q_0)|^2+|x-y|^2\\
    &\leq l(Q_0)^2 (k_1^2+2k_1+2)
\end{align}
On the other hand, we have that:
\begin{align}
    |s-\tau|&=|s-t-k_2l(Q_0)^2|\\
    &=l(Q_0)^2\left|\frac{s-t}{l(Q_0)^2}-k_2\right|\\
    &=l(Q_0)^2|\left(k_2-\frac{|s-t|}{l(Q_0)^2}\right)\\
    &\geq l(Q_0)^2(k_2-1)
\end{align}
and the choice of $a$ needed to apply the lemma becomes apparent. Doing so lets us observe that by both the Harnack inequality \ref{Harnack} and Bourgain estimate \ref{bourg}: 
\begin{equation}\w_L^{\X^+_{Q_i^l}}(Q_i^l\cap(\Tilde{O}_{j-1}\setminus O_j))\leq \frac{C_\eta}{\w_L^{\X^+_{Q_i^l}}(Q_i^l)} \w_L^{\X^+_{Q_0}}(Q_i^l\cap(\Tilde{O}_{j-1}\setminus O_j))\leq \frac{C_\eta}{\w_L^{\X^+_{Q_i^l}}(Q_i^l)} \w_L^{\X^+_{Q_0}}(Q_i^l\cap O_{j-1})\end{equation}
Recalling that the $O_j$ are a good $\epsilon_0$ cover with respect to $\w_L^{\X^+_{Q_0}}$ we make use of the previous estimates and the Harnack inequality to bound $u(\X_{\Tilde{Q}_i^l}^+)$ from above:
\begin{align}
    u(\X_{\Tilde{Q}_i^l}^+)&\leq \w_L^{\X_{\Tilde{Q}_i^l}^+}(\partial\Omega\setminus Q_i^l)+\w_L^{\X_{\Tilde{Q}_i^l}^+}(S\cap Q_i^l)\\
    &\leq C\eta^\gamma + \sum_{j=l+2}^k\w_L^{\X_{\Tilde{Q}_i^l}^+}(Q_i^l\cap(\Tilde{O}_{j-1}\setminus O_j))+\w_L^{\X_{\Tilde{Q}_i^l}^+}(Q_i^l\cap(\Tilde{O}_l\setminus O_{l+1}))\\
    &\leq C\eta^\gamma+\frac{C_\eta}{\w_L^{\X_{\Tilde{Q}_i^l}^+}(Q_i^l)}\sum_{j=l+2}^k\w_L^{\X_{Q_0}}(Q_i^l\cap O_{j-1})+\w_L^{\X_{\Tilde{Q}_i^l}^+}(\tilde{Q}_i^l)\\
    &\leq C\eta^\gamma+C_\eta\epsilon_0+1-c_0 
\end{align}
so up to selection of $\eta,\,\epsilon_0$ we can conclude $u(\X_{\Tilde{Q}_i^l}^+)\leq 1-\frac{3}{4}c_0$. On the other hand we observe:
\begin{align}
    u(\X_{\Tilde{Q}_i^l}^+)&=\w_L^{\X_{\Tilde{Q}_i^l}^+}(S)\geq \w_L^{\X_{\Tilde{Q}_i^l}^+}(Q_i^l\cap(\Tilde{O}_l\setminus O_{l-1}))\\
    &=\w_L^{\X_{\Tilde{Q}_i^l}^+}(\Tilde{Q}_i^l\setminus O_{l+1})\\
    &=\w_L^{\X_{\Tilde{Q}_i^l}^+}(\Tilde{Q}_i^l)-\w_L^{\X_{\Tilde{Q}_i^l}^+}(\Tilde{Q}_i^l\cap O_{l+1})\\
    &\geq c_0-\w_L^{\X_{\Tilde{Q}_i^l}^+}(\Tilde{Q}_i^l\cap O_{l+1})\\
    &\geq c_0-C_\eta \epsilon_0
\end{align}
and again, up to choosing $\eta,\,\epsilon_0$ to be sufficiently small, we can conclude that $u(\X_{\Tilde{Q}_i^l}^+)\geq \frac{3}{4}c_0$.\\
\vspace{1em}
\\
Now we estimate $u(\X_{\Tilde{P}_i^l}^+)$ from above in the case that $P_i^l\cap \Tilde{Q}_i^l = \emptyset$. In service of this we first seek to show that $\w_L^{\X_{\Tilde{P}_i^l}^+}(Q_i^l(\Tilde{O}_{j-1}\setminus O_j))\leq \w_L^{\X_{\Tilde{Q}_i^l}^+}(Q_i^l(\Tilde{O}_{j-1}\setminus O_j))$. But observing that $l(\Tilde{P}_i^l)=\eta^2l(Q_i^l)=:\eta^2l$ and $\delta(X^+_{\Tilde{P}_i^l})\geq \eta^2l\sqrt{\frac{k_1^2}{M^2+1}+k_2}$ we can readily employ the Harnack inequality by letting $a=\eta^2\frac{1}{\sqrt{2(k_1^2+k_2)}-1}$. Starting with the same estimation as before we have:
\begin{align}  
(\X^+_{\tilde{P}_i^l})&\leq C\eta^\gamma +\sum_{j=l+2}^k\w_L^{\X^+_{\tilde{P}_i^l}}(
P_i^l\cap(\Tilde{O}_{j-1}\setminus O_j))+\w_L^{\X^+_{\tilde{P}_i^l}}(
P_i^l\cap(\Tilde{O}_{l}\setminus O_{l+1}))\\
&\leq C\eta^\gamma+C_{\eta^2} \epsilon_0+\w_L^{\X^+_{\tilde{P}_i^l}}(
P_i^l\cap(\Tilde{O}_{l}\setminus O_{l+1}))\\
&=C\eta^\gamma+C_{\eta^2}\epsilon_0
\end{align}
Where we conclude that the third term on the second line is zero as $P_i^l\cap(\Tilde{O}_{l}\setminus O_{l+1})=(P_i^l\cap\Tilde{Q}_i^l)\setminus O_{l+1}=\emptyset$. Now, for sufficiently small $\eta$, we have that $u(\X_{\Tilde{P}_i^l}^+)<\frac{1}{4}c_0$ thus:
\begin{align}
    \left|u(\X_{\Tilde{Q}_i^l}^+)-u(\X_{\Tilde{P}_i^l}^+)\right|&=u(\X_{\Tilde{Q}_i^l}^+)-u(\X_{\Tilde{P}_i^l}^+)\geq \frac{1}{2}c_0
\end{align}
\\
\vspace{1em}
\\
To cover the case where $P_i^l\cap \Tilde{Q}_i^l\neq \emptyset$, we estimate $u(X^+_{\Tilde{P}_i^l})$ from below. To do this, we first note that $l(P_i^l)=l(\Tilde{Q}_i^l)$. This, along with the fact that they're both dyadic cubes and non-disjoint, means that $P_i^l=Q_i^l$. Moreover, this means that $P_i^l\cap\Tilde{O}_l=\Tilde{Q}_i^l\cap \Tilde{O}_l=\Tilde{Q}_i^l=P_i^l$, thus:
\begin{align*}
    u(\X_{\Tilde{P}_i^l}^+)&=\w_L^{\X_{\Tilde{P}_i^l}^+}(S)\\
    &\geq \w_L^{\X_{\Tilde{P}_i^l}^+}(P_i^l\cap(\Tilde{O}_l\setminus O_{l+1}))\\
    &=\w_L^{\X_{\Tilde{P}_i^l}^+}(P_i^l\setminus O_{l+1})\\
    &=\w_L^{\X_{\Tilde{P}_i^l}^+}(P_i^l)-\w_L^{\X_{\Tilde{P}_i^l}^+}(P_i^l\cap O_{l+1})\\
    &\geq1-C\eta^\gamma-\w_L^{\X_{\Tilde{P}_i^l}^+}(\Tilde{Q}_i^l\cap O_{l+1})\\
    &\geq1-C\eta^\gamma-\w_L^{\X_{\Tilde{Q}_i^l}^+}(\Tilde{Q}_i^l\cap O_{l+1})\\
    &\geq 1-C\eta^\gamma-C_\eta\epsilon_0
\end{align*}
Note the inequality in the penultimate line was the Harnack inequality, passing from $\X_{\Tilde{P}_i^l}^+$ to $X^+_{\Tilde{Q}_i^l}$, which is valid in this case as $\Tilde{P}_i^l\subset  P_i^l=\Tilde{Q}_i^l$. Now, for sufficiently small $\eta$ and $\epsilon_0$ we have that $u(\X_{\Tilde{P}_i^l}^+)\geq 1-\frac{1}{4}c_0$ thus:
\begin{align}
    \left|u(\X_{\Tilde{Q}_i^l}^+)-u(\X_{\Tilde{P}_i^l}^+)\right|&=u(\X_{\Tilde{P}_i^l}^+)-u(\X_{\Tilde{Q}_i^l}^+)\geq \frac{1}{2}c_0
\end{align}   
\end{proof}
With this proof in hand we now have the tools necessary to complete the proof of lemma \ref{oncube} and, in turn, complete the proof of \ref{mainthm}.
\begin{proof}[Proof of Lemma \ref{oncube}] Take $\beta$, $Q_0$, and $F$ as given and leave $\alpha$ to be specified, but for now mandate that $0<\alpha<\alpha_0$ as prescribed by Lemma 4.7. Now let $\ulr{\X}=\w_L^{\X}(S)$ as given by 4.7. Then we have that:
\begin{align}
    \frac{-\ln(\alpha)}{C_\eta^2}\sigma(F)&\leq \int_FS_{Q_0} \ulr{\X}^2d\sigma\left(\X\right)\\
    &\leq \int_{Q_0}\left(\int_{\Gamma_{Q_0}(\tin{\X})}|\nabla u(\Y)|^2\delta(\Y)^{-n}d\Y\right)d\sigma \\
    &=\int_{B^\sigma(\tin{\X},l(Q_0))}|\nabla \ulr{\Y}|^2\delta(\Y)^{-n}\left(\int_{Q_0}\mathbf{1}_{\Gamma_{Q_0}(\tin{\X})}(\Y)d\sigma\left(\X\right)\right)d\Y
\end{align}
Now, given that $\Gamma_{Q_0} $ is a parabolic cone and $Y\in B(\X,l(Q_0))$, there exists a constant depending only on the structure of $\Omega$ for which $\{\X:\,\Gamma_{Q_0}(\X) \ni \Y\}\subset B(\hat{\Y},C\delta(\Y))\cap\partial\Omega$ where $\hat{\Y}\in\partial\Omega$ is the point below $\Y$ on $\pom$. From this, we see that for each fixed $\Y\in B(\X,l(Q_0)$, that:
\begin{equation}
    \int_{Q_0}\mathbf{1}_{\Gamma_{Q_0}(\tin{\X})}(\Y)d\sigma \leq\sigma(B^\sigma(\hat{\Y},C\delta(\Y))\cap\partial\Omega)\approx C\delta(\Y)^{n+1}
\end{equation}

Now, as a bounded solution in $\Omega$, we can apply the Carleson estimate to $u$ and see that:
\begin{align}\frac{-\ln(\alpha)}{C_\eta^2}\sigma(F)&\leq \int_{B^\sigma(\X,l(Q_0))\cap\Omega}|\nabla \ulr{\Y}|^2\delta(\Y)^{-n}\left(\int_{Q_0}\mathbf{1}_{\Gamma_{Q_0}(\tin{\X})}(\Y)d\sigma\left(\X\right)\right)d\Y\\
&\leq \int_{B^\sigma(\X,l(Q_0))\cap\Omega}|\nabla \ulr{\Y}|^2\delta(\Y)^{-n}C\delta(\Y)^{n+1}d\Y\\
&=\int_{B^\sigma({\X},l(Q_0))}|\nabla \ulr{\Y}|^2\delta(\Y)d\Y\\
&<C\sigma(B^\sigma(\X,l(Q_0))\cap\Omega)\\
&\approx C\sigma(Q_0)
\end{align}
After choosing a sufficiently small $\alpha$, this leaves us with $\displaystyle\frac{\sigma(F)}{\sigma(Q)}\leq \beta$, as desired.

\end{proof}
\section{$A_\infty$ Implies CME}
In this section we prove the direction:
\begin{theorem}[$A_\infty$ implies CME]\label{mainthm2}
    Let $\om\in\ree$ be a a Lipschitz graph domain, and let $L$ be a uniformly-parabolic operator.
    If for every $\beta\in(0,1)$ there will exist an $\alpha\in(0,1)$ such that for every surface cube $\Delta_0$, surface subcube $\Delta\subset\Delta_0$ and Borel set $F\subset \Delta$ the implication \eqref{4a} holds, then every bounded solution $u$ to the Dirichlet problem corresponding to $L$ on $\om$ satisfies \eqref{CME}.
\end{theorem}
Note that \eqref{CME} only mandates that we find a uniform bound given any $\X\in\pom$ and any $I$ with center $\X$. So for some $\X_0\in\pom$, fix a cube $I_0=I(X_0,R)$ centered at $\X_0$ and let $\Delta_0=I_0\cap\pom$. For arbitrary $r\in\re^+$ we define $\Delta_r=I(\X_0,r)\cap\pom$, noting that $\Delta_0\equiv\Delta_R$. Now . In what follows, for an open set $S$ we will denote the whitney decomposition of $S$ as $\W(S)$. Lastly, let $\w=\w_L^{\X_{\Delta_{C_0R}}^+}$ and $\Tilde{G}(\Y)=G(\X_{\Delta_{C_0R}}^+,\Y)$ for some sufficiently large $C_0>0$ to be chosen later.\\
\par Now, if we seek to show that:
$$\int_{I_0\cap\om}|\nabla u |^2\delta(\X) d\X\leq C_0\sigma(\Delta_0)$$
it suffices to show that: $$\int_{\om'}|\nabla u |^2\delta(\X) d\X\leq C_0\sigma(\Delta_0)$$
for some $\om'\supseteq I_0\cap\om$. In order to define this region we first define a subset of $\W(\om)$ relative to some $Q\in\dd(\pom)$: $$W_Q:=\left\{I\in \W(\om):\,\frac{1}{k}\leq \frac{\delta(I,Q)}{l(Q)}\leq k,\,l(I)\approx l(Q)\right\}$$
where $k$ is a universally fixed constant. From this subset, we define the following Whitney regions and associated quantity:
\begin{align*}
    U_Q&:=\bigcup_{W_Q^k}I^\sigma\\
    U_Q^*&:=\bigcup_{W_Q^k}I^{2\sigma}\\
    \alpha_Q&:=\int_{U_Q}|\nabla u |^2\delta d\X
\end{align*}
noting here that $\sigma$, though still a constant depending only on structural quantities, may differ from the occurences in the previous sections. Using these Whitney regions, we let $\dd(\Delta_0)$ be the collection of cubes $Q'\in\dd(\pom)$ such that $Q\subset\Delta_0$ and define the set $\om'$ to be:
$$\om':=\bigcup_{Q'\in\dd(\Delta_0)}U_{Q'}$$
and note that because each $U_{Q'}$ is a Whitney region, we have that:
$$\int_{\om'}|\nabla u |^2\delta(\X) d\X\approx\sum_{Q'\in\dd(\Delta_0)}\alpha_{Q'}$$
From here, we further refine the quantity to bound by the following discretized Nirenberg-type lemma:
\begin{lemma}[Discretized John-Nirenberg Lemma]\label{niren} Suppose there exist universal $\eta,c_0$ such that for all $Q\in\dd(\pom)$ there exists a family of subcubes $\F_Q:\{Q_j\}$ with $Q_j\subset Q\, \forall j$ such that:
\begin{align*}
    \sigma(Q\setminus\cup_{F_Q}Q_j)&\geq \eta\sigma(Q)\\
    \sum_{Q'\in S(Q)}\alpha_{Q'}&\leq c_0\sigma(Q)
\end{align*}
    where $S_Q$ is the family of cubes $S_Q:=\{Q'\dd(Q):Q'\not\subset Q_j,\,\forall Q_j\in\F_Q\}$. Then:
    $$\sum_{Q'\in\dd(Q)}\alpha_{Q'}\leq\frac{c_0}{\eta}\sigma(Q)$$
\end{lemma}
We make use of this lemma by defining $\F_Q$ and $S_Q$ relative to a closed subset $F$ of $\Delta_0$ constructed in the following lemma from \cite{bortz2024}:
\begin{lemma}\label{lemF}
    Given a surface cube $\Delta_0$, an $\epsilon>0$, and $\w\in A_\infty(\sigma)$, there exists a constant $M\geq 1$ depending only on structural constants and $\epsilon$, and a closed set $F$ such that: \begin{equation}F\subset\tilde{F}:=\left\{\X\in\Delta:\frac{1}{M}\leq \frac{\w(\Delta(\X,r))}{\sigma(\Delta(\X,r))}\leq M,\,\forall r \in(0,5Ml(\Delta_0))\right\}\end{equation}\\
    \begin{equation}
        \sigma(\Delta_0\setminus F)\leq \epsilon \sigma(\Delta_0)
    \end{equation}
\end{lemma}
With this lemma, we let $\F:=\F_{\Delta_0}$ be the maximal subcubes of $\Delta_{CR}$ that do not meet this $F\subset\Delta_0$ found in the lemma above. Correspondingly, we let $S:=S_{\Delta_0}$ be the set $S_{\Delta_0}:=\{Q\subset \Delta_0: Q\not\subset Q_j,\,\forall Q_j\in\F_{\Delta_0}\}$. Note quickly here that $\F_{\Delta_0}$ is exactly the family need to make use of \ref{niren} with $\eta=1-\epsilon$:
\begin{align}
    \sigma(\Delta_0\setminus \cup_{Q'\in\F} Q')&\geq\sigma(\Delta_0\cap F)\\
    &=\sigma(\Delta_0)-\sigma(\Delta_0\setminus F)\\
    &=(1-\epsilon)\sigma(\Delta_0)
\end{align}
Now we define the following domain: \begin{equation}\om_F:=\bigcup_{Q'\in S}U_{Q'}\end{equation}
taking note that by the bounded overlap property of Whitney regions, and \ref{niren}, we have that: \begin{equation}\label{imp}\sum_{Q'\in S}\alpha_{Q'}\approx\int_{\om_F}|\nabla u |^2\delta(\X) d\X\bestrel{<}\sigma(\Delta_0)\quad\Rightarrow\quad \sum_{Q'\subset \Delta_0}\alpha_{Q'}\approx \int_{\om'}|\nabla u |^2\delta(\X) d\X\bestrel{<}\sigma(\Delta_0)\end{equation}
This means that we have reduced matters to bounding the integral on the lefthand side of the implication in $\eqref{imp}$. In order to do this, we make the following observation:
\begin{lemma}\label{gd}
    $\X\in\om_{F}\Rightarrow \Tilde{G}(\X) \approx \delta(\X) $
\end{lemma}
The proof of this relies on the following lemma:
\begin{lemma}\label{gds} For all $\Y\in I_0\cap\Omega$, and $\w=\w_L^{\X^+_{\Delta_{C_0R}}}$ with $C_0$ sufficiently large we have that: 
    $$\frac{\Tilde{G}(\Y)}{\delta(\Y)}\approx\frac{\w(\Delta_{\Y})}{\sigma(\Delta_{\Y})}$$
\end{lemma}
\begin{proof}[Proof of \ref{gds}]
    First, we show $\bestrel{<}$. To do this fix $\Y\in Q_0$ and center a cube $Q\ni \Y$ such that $l(Q)=\gamma\delta(\Y)$ for some $\gamma<<1$ to be chosen to be as small as desired later. Now, construct domain $\om_{\tin{\Y}}$ to be the relative complement of $Q$ in $\om$ intersected with the domain $\{t\geq s\}$ and note that $\partial\om_{\Y}$ is composed of the parts $(\{t=s\}\setminus Q)\cap\om$ and $\partial Q\cap\{t\geq s\}=:\partial Q^+$. Further note that $\X_{\Delta_C}^+\in\om_Y$ for sufficiently large $C$, so we proceed by the maximum principle. Noting that $G(\X,\Y)\approx (t-s)^{\frac{-n}{2}}e^{-(\frac{|\X-\Y|}{c(t-s)})^2}$, we trivially have that $G(\X,\Y)\equiv 0$ for $\X\in(\{t=s\}\setminus Q)\cap\om$. For $\X\in\partial Q^+$, first note that by the Bourgain estimate and doubling, we have that $\w_L^{\X}(\Delta_Y)\sim 1$ and as $$G(\X,\Y)\bestrel{<}||\X-\Y||^{-n}\approx \sigma(\Delta_{\Y})^{-1}\delta(\Y)$$Thus $$\frac{\sigma(\Delta_{\Y})}{\delta(\Y)}G(\X,\Y)\bestrel{<}1\approx \w_L^{\X}(\Delta_{\Y})$$ 
    Now, we show $\bestrel{>}$. To do this, first let $J(\Y)$ be a parabolic rectangular prism centered at $\Y$ for which $\Delta_{\Y}\subset\subset J(\Y)\cap\pom$ and has side lengths comparable to $\delta(\Y)$, and construct $\Phi\in C_0^\infty(\frac{3}{2}J(\Y))$ for which $\Phi\equiv 1$ on $J(\Y)$, $|\nabla\Phi|\bestrel{<}\delta(\Y)^{-1}$, and $|\partial_t\Phi|\bestrel{<}\delta(\Y)^{-2}$. Now, noting that $\Phi(\X_{\Delta_C})=0$, Riesz' Formula yields that:
    \begin{align*}
        \w(\Delta_{\Y})&\leq\int\Phi d\w \\
        &=-\int A \nabla\tilde{G} \cdot\nabla\Phi +\tilde{G} \X \partial_t\Phi d\X\\
        &\bestrel{<}\int |\nabla\tilde{G} ||\nabla\Phi |+|\tilde{G} ||d_t\Phi d\X\\
        &\bestrel{<}\delta(\Y)^{-1}\int_{\frac{3}{2}J(\tin{\Y})}|\nabla\tilde{G} |d\X+\delta(\Y)^{-2}\int_{\frac{3}{2}J(\tin{\Y})}|\tilde{G} |d\X\\
        &\bestrel{<}\delta(\Y)^{-2}\int_{2J(\tin{\Y})}|\tilde{G} |d\X\\
        &\bestrel{<}\delta(\Y)^{-1}\tilde{G}(\Y)\sigma(\Delta_{\Y})
    \end{align*}
\end{proof}
With this in hand, we can now complete the proof of \ref{gd}
\begin{proof}[Proof of \ref{gd}]
    Given \ref{gds} and the definition of $F$, fixing some $\X\in\om_{F}$, it suffices to show that for some $\Y\in F$ and $r>0$ $$\frac{\w(\Delta_{\X})}{\sigma(\Delta_{\X})}\approx\frac{\w(\Delta(\Y,r))}{\sigma(\Delta(\Y,r))}$$
    First, take $\X\in J_{\X}\in W_Q$ for some $Q\in S$. Note this means $\delta(\X) \approx l(Q)$. By definition, we also have that $Q\cap F\neq \emptyset$, so take some $\Y_{F,Q}\in Q\cap F$ and construct $\Delta_{F,Q}:=\Delta(\Y_{F,Q},r_Q)$ with $r_Q$ chosen so that $\forall \X\in U_Q$ we have $\Delta_X\subset\Delta_{F,Q}$ and $r_Q\bestrel{<}R$ so that by \ref{lemF}, $w(\Delta_{F,Q})\approx \sigma(\Delta_{F,Q})$. This is possible; note that for $\X'\in\Delta_{\X}$ we have:
    \begin{align*}
        \delta(\Y_{F,Q},\X')&\leq \diam(Q)+\delta(Q,\Delta_{\X})+\diam(\Delta_{\X})\\
        &\leq \diam(Q)+\dist(Q,\hat{\X})+r_{\X}+\diam(\Delta_{\X})\\
        &\bestrel{<}l(Q)+\delta(Q,\X)+r_{\X}+r_{\X}\\
        &\bestrel{<}l(Q)+\delta(Q,\X) \\
        &\bestrel{<}l(Q)<R
    \end{align*}
    Now, by doubling of $\w,\sigma$ on $\Delta_0$ we have: $$1\approx\frac{\w(\Delta_{F,Q})}{\sigma(\Delta_{F,Q})}\approx\frac{\w(\Delta_{\X})}{\sigma(\Delta_{\X})}\approx \frac{\Tilde{G}(\X) }{\delta(\X) }$$
\end{proof}
\begin{proof}[Proof of \ref{mainthm2}]
    In what follows, we will seek to integrate by parts and will need to stay quantitatively far from the boundary to facilitate this. Define the following quantities for each $N\in\NN$: \begin{align}
    S_N&:=\{Q\in S: l(Q)>2^{-N}\}\\
    \om_{F,N}&:=\bigcup_{Q'\in S_N}U_{Q'}
\end{align}
Now, using standard methods from \cite{Stein} we let \begin{equation}
    \W_N:=\{J\in\W:J\in \W_{Q}, Q\in S_N\}
\end{equation}
and form a partition of unity on $\om_{F,N}$ constructed of functions $\eta_J\in C_0^\infty(J^{2\sigma})$ for each $J\in \W_N$ such that: \begin{equation}\label{eta}\eta:=\sum_{J\in \W_p}\eta_J \bestrel{<}1, \quad\forall\X\in\om_{F,N}\end{equation}
and for each $J\in \W_N$:
\begin{equation}\label{jdir}
    l(J)|\nabla\eta_J|+{l}(J)^2|\partial_t\eta_J|\bestrel{<}1
\end{equation}
Defining the domain \begin{equation}
    \om_{F,N}^*:=\bigcup_{Q'\in S_N}U_{Q'}^*
\end{equation}
we see that: \begin{equation}
    \mathbf{1}_{\om_{F,N}}\leq \eta\leq \mathbf{1}_{\om_{F,N}^*}
\end{equation}
meaning by the monotone convergence theorem we have further reduced matters to showing that:
\begin{equation}\label{inteta}
    \int_{\om_{F,N}}|\nabla u |^2 \delta(\X) d\X\leq \int_{\ree}\eta |\nabla u |^2 \delta(\X) d\X\leq C'\sigma(\Delta_0)
\end{equation}
\newcommand{\tg}{\Tilde{G}}
\newcommand{\reint}{\int_{\ree}}
\newcommand{\del}{\nabla}
\newcommand{\dt}{\partial_t}
We will do this by reducing the inequality on the right hand side of \ref{inteta} to a Carleson estimate on $\eta$ itself. We begin by first noting that: 
\begin{align*}
        \int_{\ree}\eta  |\nabla u |^2 \delta(\X) d\X&\bestrel{<}\int_{\ree}\eta  |\nabla u||A\nabla u | \delta(\X) d\X\\
    &=\frac{-1}{2}\int_{\ree}\eta L(u^2)\delta(\X) d\X\\
    &\approx-\int_{\ree}\eta \left[\nabla\cdot A\nabla u^2 -\partial_t u^2 \right]\tg d\X\\
    &\approx\reint A\nabla u^2 \cdot\nabla(\eta \tg )-u^2 \partial_t(\eta \tg )d\X\\
    &=\reint -u^2 \eta (\del\cdot A^t\del\tg +\partial_t\tg )+\tg (A\del u^2 \cdot\del\eta )-u^2 (A^t\del\tg \cdot\del \eta +\tg \dt\eta )d\X\\
    &=\reint \tg (A\del u^2 \cdot\del\eta )-u^2 (A^t\del\tg \cdot\del \eta +\tg \dt\eta )d\X
\end{align*}
as, by definition, $\tg$ is an adjoint solution to $L$. Continuing, we observe that:
\begin{align*}
    \int_{\ree}\eta  |\nabla u |^2 \delta(\X) d\X&\bestrel{<}\reint \tg (A\del u^2 \cdot\del\eta )-u^2 (A^t\del\tg \cdot\del \eta +\tg \dt\eta )d\X\\
    &\leq \reint2\tg u (A\del u\cdot\del\eta)+u^2(|A^t\del\tg||\del\eta|+|\tg||\dt\eta|)d\X\\
    &\bestrel{<}\reint |u||\del u||\tg||\del\eta|+u^2(|\del\tg|\del\eta|+||\tg||\dt\eta|)d\X\end{align*}
Now we take note that both $|\del\eta|$ and $|\dt\eta|$ can only be nonzero on the dilated Whitney cubes $J^{2\sigma}\subset\om^\star_{F,N}$ such that $J^{2\sigma}\cap(\om_{F,N}^*\setminus\om_{F,N})\neq\emptyset$. We label the collection of such cubes $\B_N$. This means we now need to show:
\begin{equation}
   \sum_{\B_N}\reint |u||\del u||\tg||\del\eta_J|+u^2(|\del\tg||\del\eta_J|+|\tg||\dt\eta|)d\X\,\bestrel{<}\,\sigma(\Delta_0)
\end{equation}
To do this, we show the following string of inequalities:
\begin{equation}\label{ecbd}
     \sum_{\B_N}\reint |u||\del u||\tg||\del\eta_J|+u^2(|\del\tg||\del\eta_J|+|\tg||\dt\eta_J|)d\X\,\bestrel{<}\,\sum_{\B_N}l(J)^{n+1}\,\bestrel{<}\,\sigma(\Delta_0)
\end{equation}
For the first inequality, we show this for each term in the integrand via the properties of Whitney cubes, the Cacciopoli inequality, and Hölder's inequality:
\newcommand{\js}{J^{2\sigma}}
\newcommand{\rjs}{\int_{J^{2\sigma}}}
\newcommand{\lrh}[1]{\left(#1\right)^{\frac{1}{2}}}
\begin{align*}
    \rjs u^2|\del\tg||\del\eta_J|d\X&\bestrel{<}_u\rjs|\del\tg||\del\eta_J|d\X
    \\
&\leq \left(\rjs|\del\tg|^2d\X\right)^{\frac{1}{2}}\left(\rjs|\del\eta_J|^2d\X\right)^{\frac{1}{2}}\\ 
    &\bestrel{<}\,\left(\frac{1}{l(J)^2}\int_{2\js}|\tg|^2d\X\right)^{\frac{1}{2}}\left(\rjs\delta(\X)^{-2}d\X\right)^{\frac{1}{2}}\\
    &\bestrel{<}\,\left(\int_{2\js}d\X\right)^\frac{1}{2} l(J)^{-1}\left(\rjs d\X\right)^\frac{1}{2}\\
   &\bestrel{<}\, l(J)^{-1}|J|\\
    &\bestrel{<}\, l(J)^{n+1}\\
    \\
    \rjs u^2|\tg||\dt\eta_J|d\X&\leq \rjs |\tg||\dt\eta_J|d\X\\
    &\bestrel{<}\, \rjs l(J)^{-1} d\X\\
    &\bestrel{<}\, l(J)^{n+1}\\
    \\
    \rjs |u||\del u||\tg||\del\eta_J|d\X&\leq \lrh{\rjs|\tg|^2|\del\eta_J|^2|u|^2d\X}\lrh{\rjs|\del u|^2d\X}\\
    &\bestrel{<}\,\lrh{\rjs \delta(\X)^2|\del\eta_J|^2d\X}\lrh{l(J)^{-2}\int_{2\js} u^2 d\X }\\
    &\bestrel{<}\,\lrh{\rjs l(J)^2|\del\eta_J|^2d\X}l(J)^{-1}|J|^{\frac{1}{2}}\\
    &\bestrel{<}\,l(J)^{-1}|J|\\
    &\bestrel{<}\,l(J)^{n+1}
\end{align*}
\newcommand{\teta}{\Tilde{\eta}}
\newcommand{\omstar}{\om_{\star\star,N}}
\newcommand{\omt}{\om_{\star\star\star}}
\newcommand{\beg}{\,\bestrel{>}\,}
\newcommand{\bel}{\,\bestrel{<}\,}
This takes care of the first inequality of \ref{ecbd}. For the second, we look to make use of a similar estimate from \cite{bortz2024}. In that paper, they define a similar cutoff function $\teta$ for a corresponding domain $\omstar$ defined as $\omstar:=\{\Y=(y_0,y,s)\in\Delta_{CR}: y_0\geq \Psi(y,s)+2^{-1}c_1 h(y,s)+2^{-N}\}$ where $h(y,s)$ is a regularized distance function to the projection $\pi(F)$ of the set $F$ onto the plane $0\times\re^{n-1}\times\re$ with the property that $h(y,s)\approx \delta(\pi(F),(y,s))$. This $\teta$ is constructed to have compact support in the set $\omt:=\{\Y=(y_0,y,s)\in\Delta_{CR}: y_0\geq \Psi(y,s)+2^{-2}c_1 h(y,s)$. Define $\Tilde{\B}_N$ to be the Whitney cubes that meet both $\omstar$ and $\omt\setminus \omstar$. It is then shown in \cite{bortz2024} that:
$$\sum_{\Tilde{\B}_N}l(J)^{n+1}\,\bestrel{<}\,R^{n+1}\,\bestrel{<}\,\sigma(\Delta_0)$$
This result can be extended to our set $\B_N$ if it can be shown that $\om_{F,N}\subseteq \omstar$. To show this, it is sufficient to show that for every $\X=(x_0,x,t)\in\om_{F,N}$ that $|x_0-\Psi(x,t)|\,\bestrel{>}\, h(x,t)+2^{-N}$. First note that by nature of a Lipschitz graph domain, it suffices to show that $\delta(\X)\,\bestrel{>}\, h(x,t)+2^{-N}$. Trivially, we have that $\delta(\X)\geq 2^{-N}$. To show that $\delta(\X)\beg h(x,t)$, first note that $\X$ is contained in some $U_Q$ for some $Q$ that meets $F$. This means there exists some point $\textbf{Z}=(z_0,z,\zeta)\in Q\cap F$ for which: $h(x,t)\bel \delta(\X,\Z)\bel\delta((x,t),(z,\zeta))\bel l(Q)\bel \delta(\X)$. 
\end{proof}
\bibliographystyle{plain}
\bibliography{refs}

\end{document}